\documentclass{amsart}
\usepackage{graphicx, amsmath, amssymb,amsthm, geometry, caption}
\usepackage{graphicx}
\usepackage{amsrefs}

\newtheorem{Thm}{Theorem}
\newtheorem{lem}{Lemma}
\newtheorem{prop}{Proposition}
\newtheorem{remark}{Remark}

\newtheorem{defn}{Definition}
\newtheorem{coro}{Corollary}

\def\*#1{\mathbf{#1}} %define \*something = \mathbf{something}
\date{}

\begin{document}
\title{A surface area formula for compact hypersurfaces in $\mathbb{R}^n$}
\author{Yen-Chang Huang}
\address{National University of Tainan, Tainan city, Taiwan R.O.C.}
\curraddr{}
\email{ychuang@mail.nutn.edu.tw}
\thanks{This work was funded by National Science and Technology Council (NSTC), Taiwan, with grant Number: 110-2115-M-024 -002 -MY2.}
%\thanks{}

%\date{5th, March 2023 (submission)}

\dedicatory{}

\begin{abstract}
The classical result of Cauchy's surface area formula states that the surface area of the boundary $\partial K=\Sigma$ of any $n$-dimensional convex body in the $n$-dimensional Euclidean space $\mathbb{R}^n$ can be obtained
by the average of the projected areas of $\Sigma$ along all directions in $\mathbb{S}^{n-1}$.
In this notes, we generalize the formula to the boundary of arbitrary $n$-dimensional submanifolds in $\mathbb{R}^n$ by defining a natural notion of projected areas along any direction in $\mathbb{S}^{n-1}$.
This surface area formula derived from the new concept coincides with not only the result of the Crofton's formula but that of De Jong \cite{de2013volume} by using tubular neighborhood.
We also define the projected $r$-volumes of $\Sigma$ onto any $r$-dimensional subspaces, and obtain a recursive formula for mean projected $r$-volumes of $\Sigma$.
\end{abstract}
\maketitle

\section{Introduction and main results}
Finding the volume of embedded submanifolds $K$ of dimensional $k$ in the Euclidean spaces $\mathbb{R}^n$ with boundary $\Sigma=\partial K$ is one of the interesting research topics in differential geometry.
There have been abundant results to measure the volumes and the surface areas of $K$.
For instance, the surface area of $K$ can be measured by restricting the Euclidean metrics to $\Sigma$ and integrating the $(k-1)$-volumes of the infinitesimal parallelotopes in the tangent space of $K$ \cite[Section 2.2]{aminov2001geometry};
it can also be obtained by limiting the ratio for the volume of $(n-k)$-ball and the $n$-dimensional Lebesgue measure of the tubular neighborhood of $\Sigma$ derived by De Jong \cite{de2013volume}.

We shall restrict our attention throughout the notes to the method developed by the approach of Integral Geometry \cite{sors2004integral} and Convex Geometry \cite{schneider2014convex}.
The efficiency of the method has been validated on the applications of Geometric Tomography and other scientific fields (see more details in \cite{gardner1995geometric}).
The basic assumption in this approach is that $K$ must be convex, since the original proof was established by a limiting process of convex polyhedrons inscribed in $K$.
When $k=n$, we observe that the convexity assumption for $K$ may be relaxed by introducing a natural concept for projected surface areas of $\Sigma$ onto any $r$-dimensional subspaces, $r\leq n-1$ (see Definition \ref{defn1}).
The new notion also gives an alternative proof for Crofton's formula, which states that the surface area of $\Sigma$ can be obtained by "counting" the number of intersections of all lines with $\Sigma$ (Lemma \ref{proj1} and Theorem \ref{crofton1}, below) and
it results in the generalization of the Cauchy's surface area formula for arbitrary $K$ (Theorem \ref{cauchy2}).
We also generalize the projected surface area to higher codimensional subspaces, namely, the projected $r$-volume ($1\leq r\leq n-1$) of $\Sigma$ onto $r$-subspaces (Definition \ref{codimprojected}),
and derive a recursive formula of mean projected area of $\Sigma$ (Theorem \ref{mainthm2}).

It seems that there are only a few results concerning about the relation between projected areas and surface areas for nonconvex boundaries in the literature.
Two closer concepts probably are the integral geometric measure defined by J. Favard \cite{favard1932definition} and more recently by Bouafia-Pauw \cite{bouafia2015integral}; both are based on the approach of geometric measure theory. The kinematic formula have the similar results without introducing the concept of projections, for instance, \cites{zhang1999dual, zeng2016new}. As mentioned before, De Jong \cite{de2013volume} gave a geometric definition for volume of $r$-dimensional submanifolds in $\mathbb{R}^n$ ($1\leq r\leq n-1$), and derived a $r$-volume formula by considering the ratio of $r$-volume of the tubular neighborhood and the volume of unit ball in the normal bundle. He claimed that the formula holds when the submanifold is of dimension $n$ with $C^1$-boundary. However, we will take parallel transformations of Lie groups and the method of moving frames to construct the $r$-volume forms on the $r$-subspaces to represent the $r$-projected volumes of $\Sigma$ (see the discussion in \textit{Case 1} and \textit{Case 2} of Section \ref{content2}). In Proposition \ref{samesurfacearea}, we will also prove that the surface area we derived (see Theorem \ref{cauchy2}) coincides with that in his work \cite{de2013volume}.

The Cauchy's surface area formula in $\mathbb{R}^n$ states that the surface area of a convex hypersurface $\Sigma$ in $\mathbb{R}^n$ can be represented by the average of the projected areas of $\Sigma$ along
all normal directions of the $(n-1)$-dimensional unit sphere $\mathbb{S}^{n-1}$:

\begin{Thm}[Theorem 5.5.2 \cite{klain1997}]\label{cauchyrn}
Let $K\subset \mathbb{R}^n$ be an $n$-dimensional convex body (i.e. convex set with nonempty interior) with rectifiable boundary $\partial K=\Sigma$. The surface area (or call the $(n-1)$-dimensional volume) of $K$, denoted by $\mathcal{V}_{n-1}(\Sigma)$, is given by
\begin{align}\label{cauchy1}
{\mathcal{V}}_{n-1}(\Sigma)=\frac{1}{\omega_{n-1}}\int_{\nu\in\mathbb{S}^{n-1}}\mathcal{V}_{n-1}(\Sigma| \nu^\perp) dS_\nu.
\end{align}
Here $dS_\nu$ is the surface area element at $\nu\in \mathbb{S}^{n-1}$, $\omega_{n-1}$ is the $(n-1)$-dimensional volume of the unit ball in $\mathbb{R}^{n-1}$, and
$\mathcal{V}_{n-1}(\Sigma|\nu^\perp)$ is the $(n-1)$-dimensional volume of the orthogonal projection of $\Sigma$ onto the $(n-1)$-dimensional subspace $\nu^\perp$ perpendicular to the unit outward normal $\nu$.
\end{Thm}

The orthogonal projected area $\mathcal{V}_{n-1}(\Sigma| \nu^\perp)$ of $\Sigma$ can be explicitly represented by the integral formula
\begin{align}\label{projecteda}
\mathcal{V}_{n-1}(\Sigma| \nu^\perp)=\frac{1}{2}\int_{\tilde{\nu}\in \mathbb{S}^{n-1}} |\nu \cdot \tilde{\nu }|dS_{\tilde{\nu}},
\end{align}
where $dS_{\tilde{\nu}}$ is the surface area element at $\tilde{\nu}\in \mathbb{S}^{n-1}$.
The identities \eqref{cauchy1} and \eqref{projecteda} both can be proved by inscribing a convex $m$-polyhedron in $\Sigma$, and calculating the limit of the surface area of the polyhedron when $m$ goes to infinity.
The convexity of the polyhedron ensures that at any point in $\nu^\perp$ the counting multiplicity of the projected areas of $\Sigma$ is almost everywhere two (we may assume the polyhedron contains the origin and each point in $\nu^\perp$ is projected for twice: from the "front" and the "back" of $\nu^\perp$ respectively).
See the proof of Theorem 5.5.2 in \cite{klain1997} or \cite[page 217]{sors2004integral} for the detail. However, such method does not work for nonconvex $K$ since it may not necessarily have the inscribed convex polyhedron.

Although, in general, Theorem \ref{cauchyrn} does not work for arbitrary hypersurface $\Sigma$, the formula \eqref{cauchy1} still holds in a kind of special nonconvex case: when $\Sigma=\partial K$ is obtained such that its complement $K^c=\mathbb{R}^n \setminus K$ is convex. Indeed, suppose $K$ be a nonconvex set with convex complement $K^c$, by Lemma \ref{existconvex} (see below) we will show that the closure $\overline{K^c}$ of $K^c$ is convex. Since $\Sigma$ shares the same boundary with
$K$ and $K^c$, the surface area of $\Sigma=\partial K$ can be obtained by applying formula \eqref{cauchy1} to the convex part $\partial K^c$.

As a result, if the hypersurface $\Sigma$ obtained from a nonconvex subset $K$ such that the complement $K^c$ is also nonconvex, the argument in the previous paragraph fails.
A simple example is that when $\Sigma$ is the boundary of a star domain $S$, its complement $\mathbb{R}^n\setminus S$ is again nonconvex.
The main aim of the present paper is to generalize the notion of projected areas for nonconvex boundary and derive the Cauchy's surface area formula \eqref{cauchy1} for arbitrary boundary.

Let us make some remarks about Theorem \ref{cauchyrn}. First, when $K=\mathbb{S}^{n-1}$, the following lemma shows that the surface area of the orthogonal projection of the unit sphere is independent of the choices of the projected direction.
\begin{lem}[Lemma 5.5.1 \cite{klain1997}]\label{lem2}
For any $\nu\in \mathbb{S}^{n-1}\subset \mathbb{R}^n$,
\begin{align}
\mathcal{V}_{n-1}(\mathbb{S}^{n-1}| \nu^\perp)=\frac{1}{2}\int_{\tilde{\nu}\in \mathbb{S}^{n-1}} | \nu  \cdot \tilde{\nu}|dS_{\tilde{\nu}}=\omega_{n-1}.
\end{align}
\end{lem}
This natural property for $\mathbb{S}^{n-1}$ plays the key role for the proof of Theorem \ref{cauchyrn} (see \cite{klain1997} or the proof in Theorem \ref{cauchy2}).
We also notice that the converse of Lemma \ref{lem2} is not true in general. There exists some compact hypersurfaces in $\mathbb{R}^n$ which are not the standard spheres but with constant projected areas.
For instance, Reuleaux triangles in $\mathbb{R}^2$ and bodies of constant width in $\mathbb{R}^n$ (see \cite{chakerian1966sets}) both are not round spheres, but have constant projected areas.

Secondly, it is known that the Gauss map $\nu:\Sigma=\partial K \rightarrow \mathbb{S}^{n-1}$ is bijective if and only if $K$ is convex. Moreover, if $K$ is strictly convex, then the Gauss map is a diffeomorphism.
Since the surfaces $\Sigma$ considered in Theorem \ref{cauchyrn} and Lemma \ref{lem2} are convex, the domains over which both integrations are performed are $\mathbb{S}^{n-1}$,
which can be identified (via the diffeomorphism $\nu$) with $\Sigma$.
In consequence, a natural generalization of Theorem \ref{cauchyrn} to nonconvex hypersurfaces $\Sigma$ can be considered as an integral over $\Sigma$ itself.

We recall some fundamental background for the $n$-dimensional Euclidean space $\mathbb{R}^n$.
The space $\mathbb{R}^n$ can be regarded as an $n$-dimensional Lie group with the natural left translation $L_q$ defined by $L_q p=q+p$ for any points $q,p\in \mathbb{R}^n$.
Its inverse translation and compositions are defined to be $(L_q)^{-1}r=L_{q^{-1}}r$ and $(L_q\circ L_{p^{-1}}) r=L_{qp^{-1}}r=q-p+r$, respectively, for any $r\in \mathbb{R}^n$.
For each point $p\in \mathbb{R}^n$, we may identify the whole space $\mathbb{R}^n$ with the tangent space $T_p\mathbb{R}^n$.
In this paper, we will abuse the notations for any point $q\in \mathbb{R}^n$ and the vector $\vec{Oq}$ starting from the origin $O$ and ending at $q$ if the content is clear.
Thus, for any point $p\in \mathbb{R}^n$ and any tangent vector $q\in T_0\mathbb{R}^n$, the push-forward of the left translation, $L_{p_*}q$, is the vector obtained by moving the vector $\vec{Oq}$ to the vector starting from $p$ and parallel to $\vec{Oq}$.
Similarly, $(L_p^*)(\omega(q))$ denotes the pull-back of the $r$-form $\omega$ at $q$.
We insist using the Lie group notations in $\mathbb{R}^n$ (for instance, using $L_p q$ instead of $p+q$ for vector addition) because we believe all contents in the paper can be applied to any smooth Lie group without change. In particular, recently the author also proved the Cauchy's surface area formula for domains in the Heisenberg groups, which is regarded as a strictly pseudoconvex CR manifold with Tanaka-Webster curvature vanished (see \cite{huang2021cauchy}).
We also use the notation $n(p)\in \mathbb{S}^{n-1}$ to represent the outward unit normal $n(p)$ at $p\in \mathbb{S}^{n-1}$, and $n(p)^\perp$ for the subspace perpendicular to $n(p)$.

Let us denote by $K$ an $n$-dimensional compact submanifold embedded in $\mathbb{R}^n$ ($n\geq 1$) with rectifiable boundary $\Sigma=\partial K$ in the usual topology of $\mathbb{R}^n$. Next we give the definition for projected areas of compact hypersurfaces (not necessarily convex).
\begin{defn}\label{defn1}
The (orthogonal) weighted projected area of a compact hypersurface $\Sigma$ in $\mathbb{R}^n$ along the direction $n(p)\in \mathbb{S}^{n-1}$ onto the subspace $n(p)^\perp$ is given by
\begin{align}\label{gendefn1}
\mathcal{V}_{n-1}(\Sigma|n(p)^\perp)=\frac{1}{2}\int_{q\in \Sigma} |L_{{qp^{-1}}_*}{n}(p)\cdot \tilde{n}(q)|d\Sigma_q,
\end{align}
where $L_{{qp^{-1}}_*}{n}(p)$ is the push-forward of the outward unit normal ${n}(p)$ at $p\in \mathbb{S}^{n-1}$, $\tilde{n}(q)$ is the outward unit normal at $q\in \Sigma$, and $d\Sigma_q$ is the area element of $\Sigma$ at $q\in \Sigma$.
\end{defn}

\begin{remark}\label{geometry1}
We give a geometric interpretation for \eqref{gendefn1} as follows. At each point $q\in \Sigma$, the integrand with area element at $q$,  $|L_{{qp^{-1}}_*}{n}(p)\cdot \tilde{n}(q)|d\Sigma_q$, is the projected infinitesimal area element of $d\Sigma_q$ onto the projected point of $q$ on the plane $n(p)^\perp$ and the integral becomes the projected area of $\Sigma$ onto $n(p)^\perp$ counted with multiplicity. In other words, suppose $\ell$ is a line parallel to $n(p)$ such that $\ell\cap \Sigma\neq \emptyset$ and $\ell\cap n(p)^\perp\neq \emptyset$, all points at $\ell\cap \Sigma$ are projected along $\ell$ onto one point at $\ell\cap n(p)^\perp$ with multiple times depending on the number of $\ell\cap \Sigma$. This is the reason we call the integral the weighted projected area. When $\Sigma$ is convex, the number of the points on $\ell\cap \Sigma$ is almost everywhere two, but in general the number depends on the projected direction and $\Sigma$.
\end{remark}

According to Remark \ref{geometry1}, a geometric equivalent definition of Definition \ref{defn1} will be given in the following lemma, which states that the value $\mathcal{V}_{n-1}(\Sigma|n(p)^\perp)$ is equal to the integral of the number of intersections over all lines parallel to $n$.
\begin{lem}\label{proj1}
Given a compact hypersurface $\Sigma$ in $\mathbb{R}^n$ and a unit vector $n(p)\in \mathbb{S}^{n-1}$. Suppose $\ell$ is any line parallel to $n(p)$, and intersects with the orthogonal complement $n^\perp(p)$ at the point $u$, then the projected area of $\Sigma$  onto the subspace $n(p)^\perp$ in Definition \ref{defn1} can be obtained by
\begin{align}
\mathcal{V}_{n-1}(\Sigma|n(p)^\perp)=\frac{1}{2} \int_{\ell\perp n(p)^\perp}  \#(\ell\cap \Sigma) dn^\perp_u,
\end{align}
where $\#(\ell\cap \Sigma)$ is the number of intersections of $\ell$ and $\Sigma$; $dn^\perp_u$ is the area element of $n(p)^\perp$ at $u\in n(p)^\perp\cap\ell$.
\end{lem}

\begin{proof}
For simplicity, we fix any point $p\in \mathbb{S}^{n-1}$, and write $n=n(p)$ and $n^\perp=n(p)^\perp$. Suppose the line $\ell$ parallel to $n$ intersects $\Sigma$ at the point $q$. Then $\ell$ must intersect the subspace $n^\perp$ at the unique point $u$. We may choose two orthonormal frames $\{q, e_i\}$ and $\{u, \bar{e}_j\}$, $1\leq i,j\leq  n$, at $q$ and $u$ respectively satisfying the conditions
\begin{align}\label{condi1}
\left\{
\begin{array}{ll}
\text{at the point }q: & e_1 \in T_q\Sigma^\perp, e_2, e_\alpha\in T_q\Sigma, \\
\text{at the point }u: &\bar{e}_1 \in \ell, \bar{e}_2, \bar{e}_\alpha \in T_u n^\perp, \\
e_\alpha =\bar{e}_\alpha &\text{ for } \alpha=3,\cdots, n.
\end{array}
\right.
\end{align}
According to the construction of the frames, the transition matrix $[a^{i}_{j}]$ is given by $\bar{e}_j =\Sigma_{i=1}^n e_i a_j^i$, namely, where
\begin{align*}
\left\{
\begin{array}{lll}
a^1_1&=a^2_2 \ \  =\cos\theta, \\
a^1_2&=-a^2_1=-\sin\theta, \\
a^\alpha_\beta&=\delta^\alpha_\beta, \text{ the Kronecker delta}, 3\leq \alpha, \beta \leq n,
\end{array}
\right.
\end{align*}
and $\theta$ is the angle between the planes $T_q\Sigma$ and $n^\perp$.
Suppose $\{\omega^1, \omega^2, \omega^\alpha\}$ and $\{\bar{\omega}^1, \bar{\omega}^2, \bar{\omega}^\alpha \}$ are the dual forms of $\{e_i\}$ and $\{\bar{e}_i\}$ respectively. It can be shown that $\bar{\omega}^i=\Sigma_{j=1}^n b_j^i\omega^j$ where the matrix
$[b_i^j]=[a_i^j]^{-1}$. Therefore, by identifying the point $q\in\Sigma$ and $u\in n^\perp$, the relation between the area elements $d\Sigma_q$ and $dn^\perp_u$ at $q$ and $u$ respectively can be obtained as follows:
\begin{align}\label{transition1}
dn^\perp_u&= \bar{\omega}^2 \wedge \bar{\omega}^3\wedge \cdots \wedge \bar{\omega}^n \\
&=(\sin\theta \omega^1+ \cos\theta \omega^2)\wedge \omega^3 \wedge \cdots \wedge \omega^n  \nonumber \\
&=\sin\theta \omega^1 \wedge \omega^3 \wedge \cdots \wedge \omega^n + \cos\theta \omega^2\wedge \omega^3\wedge \cdots \wedge \omega^n  \nonumber \\
&=\sin\theta \omega^1 \wedge \omega^3 \wedge \cdots \wedge \omega^n + \cos\theta d\Sigma_q. \nonumber
\end{align}
When restrict on $\Sigma$, by \eqref{transition1} we have the projected formula
\begin{align}
dn^\perp_u \big|_{\Sigma} = |\cos\theta|d\Sigma_q
\end{align}
(here we have put the absolute value on the cosine to get the positive surface area). By integrating over all lines $\ell$ parallel to $n$, we have
\begin{align*}
\mathcal{V}_{n-1}(\Sigma|n^\perp)=\frac{1}{2}\int_{q\in \Sigma} |L_{{qp^{-1}}_*} n(p)\cdot \tilde{n}(q)| d\Sigma_q=\frac{1}{2}\int_{q\in \Sigma} |\cos\theta| d\Sigma_q=\frac{1}{2}\int_{u\in \ell\cap n^\perp} \#(\ell\cap \Sigma)dn_u^\perp
\end{align*}
and complete the proof.
\end{proof}

Recall that the Crofton's formula in $\mathbb{R}^2$ states that the perimeter of a rectifiable plane curve $\gamma$ is equal to the integral of the number $\#(\ell\cap \gamma)$ of intersections for $\gamma$ and any line $\ell$.
It has been generalized to the higher dimensional with a variety of versions (see \cite[Chapter 14]{sors2004integral} or \cite{tsukerman2017brunn}).
Notice that, in contrast to the Cauchy's surface formula, the Crofton's formula does not need the convexity assumption for curves and hypersurfaces,
and hence $\mathcal{V}_{n-1}(\Sigma|n^\perp(p))$ in Lemma \ref{proj1} seems a reasonable definition connecting both Cauchy's and Crofton's formulas.

The previous result, Lemma \ref{proj1}, will give a simpler proof for the Crofton's formula in $\mathbb{R}^n$ for $n\geq 1$. Before that, let us introduce some basic settings.
Let $\mathcal{L}$ be the set of oriented lines in $\mathbb{R}^n$.
Suppose $\ell\in \mathcal{L}$ and $\ell^\perp$ is its orthogonal complement through the origin.
$\ell$ can be uniquely determined by the following process: first, choose a line $\ell'$ parallel to $\ell$ through the origin and use the unit vector $n(p)\in \mathbb{S}^{n-1}$ to represent the direction of $\ell'$ for some $p\in \mathbb{S}^{n-1}$.
Secondly, parallel move $\ell'$ to $\ell$ at the point $u=\ell\cap \ell^\perp$.
By identifying $\ell^\perp$ with $\mathbb{R}^{n-1}$, there is a natural one-to-one correspondence between the set $\mathcal{L}$ and $\mathbb{S}^{n-1}\times \mathbb{R}^{n-1}$,
\begin{align*}
\begin{matrix}
\mathcal{L} & \longleftrightarrow & \mathbb{S}^{n-1} &\times &  \mathbb{R}^{n-1}, \\
\ell &  \longleftrightarrow & (n(p)&,& u).
\end{matrix}
\end{align*}
Thus, we may take the $(2n-2)$-form $d\ell:=dS_p\wedge d\ell^\perp_u$ as an invariant measure on $\mathcal{L}$,
where $dS_p$ is the area element at $p\in \mathbb{S}^{n-1}$ and $d\ell^\perp_u$ is the area element of $\ell^\perp$ at $u=\ell\cap\ell^\perp$. Notice that $dS_p\wedge d\ell^\perp_u$ is invariant under the rigid motions (rotations and translations) in $\mathbb{R}^n$.

\begin{Thm}[The Crofton's formula in $\mathbb{R}^n$]\label{crofton1}
Given an $n$-dimensional compact submanifold $K$ with rectifiable boundary $\partial K=\Sigma$. Let $\mathcal{L}$ be the set of all oriented lines in $\mathbb{R}^n$ and denote $\#(\ell\cap \Sigma)$ by the number of intersections of $\ell\in \mathcal{L}$ and $\Sigma$.
Then the surface area $\mathcal{V}_{n-1}(\Sigma)$ of $K$ is given by
\begin{align}
\mathcal{V}_{n-1}(\Sigma) =\frac{1}{2 \omega_{n-1} }\int_{\ell\in \mathcal{L}} \#(\ell\cap \Sigma) d\ell,
\end{align}
where $d\ell=dS_p\wedge d\ell^\perp_u$ is the invariant measure in $\mathcal{L}$ consisting of the area elements $dS_p$ at $p\in \mathbb{S}^{n-1}$ and $d\ell^\perp_u$ at $u=\ell\cap \ell^\perp$, with the orthogonal complement $\ell^\perp$ of $\ell$ through the origin.
\end{Thm}

\begin{proof}
For any line $\ell=\ell_p\in \mathcal{L}$ with the direction $n(p)\in \mathbb{S}^{n-1}$ and $\ell\cap \Sigma=q\neq \emptyset $, let $\phi_{p,q}$ be the angle between $\ell_p$ and the normal vector at $q\in \Sigma$. Notice that $d\ell_u^\perp = dn_u^\perp$ in Lemma \ref{proj1}. Using Lemma \ref{lem2}, Definition \ref{defn1}, and Lemma \ref{proj1}, a straight-forward derivation implies that
\begin{align*}
2 \omega_{n-1}\mathcal{V}_{n-1}(\Sigma)
&=2\int_{q\in \Sigma} \omega_{n-1} d\Sigma_q =\int_{q\in \Sigma}\int_{p\in \mathbb{S}^{n-1}} |\cos\phi_{p,q}|dS_p d\Sigma_q \\
&=\int_{p\in \mathbb{S}^{n-1}} \int_{q\in \Sigma} |\cos\phi_{p,q}| d\Sigma_q dS_p \\
&=\int_{p\in \mathbb{S}^{n-1}} 2 \mathcal{V}_{n-1}(\Sigma |n(p)^\perp) \ dS_p \\
&=\int_{p\in \mathbb{S}^{n-1}} \int_{\ell\perp \ell_p^\perp} \#(\ell\cap\Sigma) d\ell_u^\perp dS_q\\
&=\int_{\ell\in \mathcal{L}} \#(\ell\cap \Sigma) d\ell.
\end{align*}
\end{proof}

The following lemma is one of the motivations for the paper.

\begin{lem}\label{existconvex}
If a compact subset $K$ and its complement $K^c$ in $\mathbb{R}^n$ both are convex, then the boundary $\Sigma=\partial K$ is a hyperplane.
%Given a compact hypersurface $\Sigma$ in $\mathbb{R}^n$. If $\Sigma$ and the closure $\overline{\Sigma^c}$ of the complement of $\Sigma$ both are convex, then $\Sigma$ is a hyperplane.
\end{lem}
\begin{proof}
Since $K$ and the closure $\overline{K^c}$ of its complement $K^c$ both are closed and convex, they can be written respectively as the intersection of a family $\mathcal{F}$ and $\mathcal{G}$ of closed halfspaces, namely,
$K=\cap_{f\in\mathcal{F}}f$ and $\overline{K^c}=\cap_{g\in\mathcal{G}}g$.
We claim that all elements $f\in \mathcal{F}$ and $g\in \mathcal{G}$ are parallel. By parallel, we mean that the hyperplanes of two halfplanes are parallel.
Indeed, let us fix an element $f\in \mathcal{F}$ and assume that there exists $g\in \mathcal{G}$ such that $g$ is not parallel to $f$.
On the one hand, we have $f^c\cap g^c=(f\cup g)^c\neq \emptyset$. On the other,  $\mathbb{R}^n=K\cup \overline{K^c}\subset f\cup g$ implies that $(f\cup g)^c=\emptyset$, and we get a contradiction.
Thus, for the fixed $f$, all elements in $\mathcal{G}$ are parallel to $f$. Since $f$ can be arbitrarily chosen, we conclude that all elements in $\mathcal{F}$ are parallel to that of $\mathcal{G}$, hence $\Sigma$ is a hyperplane.
\end{proof}

Next we prove the main theorem which states that the surface area of arbitrary compact hypersurface in $\mathbb{R}^n$ is the average of the integrals of weighted projected areas over the unit sphere.
\begin{Thm} \label{cauchy2}
Let $K$ be any compact $n$-dimensional subset in $\mathbb{R}^n$ with boundary $\Sigma=\partial K$. Then its surface area is given by
\begin{align*}
\mathcal{V}_{n-1}(\Sigma)=\frac{1}{\omega_{n-1}}\int_{q\in \mathbb{S}^{n-1}} \mathcal{V}_{n-1}(\Sigma|n(q)^\perp) dS_q.
\end{align*}
\end{Thm}
\begin{proof}
By using \eqref{projecteda} and Definition \ref{defn1}, we immediately have
\begin{align*}
\int_{p\in \mathbb{S}^{n-1}}\mathcal{V}_{n-1}(\Sigma|n(p)^\perp) d S_p&=\frac{1}{2}\int_{p\in \mathbb{S}^{n-1}} \int_{q\in \Sigma}|\tilde{n}(q)\cdot {n}(p)| d\Sigma_q dS_p \\
&=\frac{1}{2}\int_{q\in \Sigma}\int_{p\in \mathbb{S}^{n-1}} |\tilde{n}(q)\cdot {n}(p)| dS_p d\Sigma_q\\
&=\int_{q\in \Sigma} \omega_{n-1} d\Sigma_q\\
&=\omega_{n-1} \mathcal{V}_{n-1}(\Sigma),
\end{align*}
and the result follows.
\end{proof}

An immediate application of Theorem \ref{cauchy2} is that a hypersurface with the smaller projected areas onto all hyperplanes has the smaller surface area; particularly, we obtain a comparison theorem of projected surface areas between two compact hypersurfaces.
Notice that, by the Alexandrov's  projection theorem \cite[page 115, Theorem 3.3.6]{gardner1995geometric}, even if two convex bodies in $\mathbb{R}^n$ have the same projected areas in all directions,
they may be completely different. In fact, there exist noncongruent convex polytopes $\Sigma_i$, $i=1,2$, with $\mathcal{V}_{n-1}(\Sigma_1 |n(p)^\perp)=\mathcal{V}_{n-1}(\Sigma_2|n(p)^\perp)$ for all $p\in \mathbb{S}^{n-1}$ (see \cite[page 121, Theorem 3.3.17]{gardner1995geometric}).
Thus, the following corollary gives a necessary condition to determine the consistency of surface areas for two hypersurfaces, but not their congruence.
\begin{coro}
Given two compact hypersurfaces $\Sigma_i$, $i=1,2$, in $\mathbb{R}^n$ and if $\mathcal{V}_{n-1}(\Sigma_1|n(p)^\perp)\leq \mathcal{V}_{n-1}(\Sigma_2 | n(p)^\perp)$ for all $p\in \mathbb{S}^{n-1}$, then $\mathcal{V}_{n-1}(\Sigma_1)\leq \mathcal{V}_{n-1}(\Sigma_2)$.
In particular, if $\mathcal{V}_{n-1}(\Sigma_1|n(p)^\perp)=\mathcal{V}_{n-1}(\Sigma_2 | n(p)^\perp)$ for all $p\in \mathbb{S}^{n-1}$, then $\mathcal{V}_{n-1}(\Sigma_1)= \mathcal{V}_{n-1}(\Sigma_2)$.
\end{coro}

Finally, we will prove that the surface area formula in Theorem \ref{cauchy2} is equivalent to that derived by De Jong in \cite{de2013volume}. More precisely, given a $k$-dimensional submanifold $M\subset \mathbb{R}^n$ and a compact subset $A\subset M$, De Jong gave a simpler method to prove the limit
\begin{align}\label{limit1}
\lim_{\epsilon \downarrow 0}\frac{\mu_n(Tub_\epsilon A)}{\beta_{n-k}\epsilon^{n-k}}
\end{align}
exits and it can be used to define the $k$-dimensional volume of $A$. Here $Tub_\epsilon A=\{p+a; p\in A, a\in N_pM, |a|<\epsilon \}$, $N_pM$ denotes the orthogonal complement of the tangent space $T_pM$ at $p$, $\mu_n(\cdot)$ is the $n$-dimensional Lebesgue measure, and $\beta_{n-k}$ is the $(n-k)$-dimensional volume of the unit ball in $\mathbb{R}^{n-k}$ (for instance, $\beta_1=2, \beta_2=\pi$, $\beta_3=\frac{4\pi}{3}$). In particular, when $A=M:=\Sigma=\partial K$ for some smooth $n$-dimensional compact submanifold $K$ in $\mathbb{R}^n$ (namely, $k=n-1$ in \eqref{limit1}), the surface area of $K$ can be obtained by
\begin{align}\label{limitsurf}
\mathcal{V}_{n-1}(\Sigma)=\lim_{\epsilon\downarrow 0}\frac{\mu_n(Tub_\epsilon \Sigma)}{2\epsilon},
\end{align}

\begin{prop}\label{samesurfacearea}
The surface area obtained in \eqref{limitsurf} for any $n$-dimensional submanifold $K$ with smooth boundary $\partial K=\Sigma$ is equal to the one obtained in Theorem \ref{cauchy2}.
\end{prop}

\begin{proof}
Since locally $\Sigma$ can be represented by a smooth defining function, we may assume that for any open subset $W\subset \Sigma$ there exists a smooth function $\phi=(\phi_1,\cdots, \phi_n):U\subset \mathbb{R}^{n-1}\rightarrow W\subset \mathbb{R}^n$ such that $W=\{\phi_n(x_1, \cdots, x_{n-1})=0; \text{ any point }(x_1, \cdots, x_{n-1})\in U\}$. Let $B=\Big(\frac{\partial \phi_i}{\partial x_j}\Big)$ for $1\leq i\leq n$ and $1\leq j\leq n-1$ be an $n\times (n-1)$-matrix.
Then the surface area element $d\Sigma_q$ at $q\in W$ satisfies
\begin{align*}
d\Sigma_q=\sqrt{\det(B^T \cdot B)} dx_1\cdots dx{_{n-1}},
\end{align*}
where $B^T$ is the matrix transpose. Thus, by Definition \ref{defn1}
\begin{align}\label{gramdeter}
\mathcal{V}_{n-1}(W|n(p)^\perp)&=\frac{1}{2}\int_{q\in W} |L_{{qp^{-1}}_*}n(p)\cdot \tilde{n}(q)|d\Sigma_q\\
&=\frac{1}{2}\int_U |\cos\theta|\sqrt{\det(B^T\cdot B)} dx_1\cdots dx_{n-1}, \nonumber
\end{align}
where $\theta$ is the angle between the unit normal vector $\tilde{n}(q)$ of $W$ and $n(p)$.
Therefore, by \eqref{gramdeter}
\begin{align}\label{localarea}
\int_{p\in \mathbb{S}^{n-1}} \mathcal{V}_{n-1}(W|n(p)^\perp)dS_{p}&=\frac{1}{2}\int_{p\in \mathbb{S}^{n-1}}\int_U |\cos\theta|\sqrt{\det(B^T\cdot B)}dx_1\cdots dx_{n-1} dS_p\\
&=\frac{1}{2}\int_U 2\omega_{n-1}\sqrt{\det(B^T\cdot B)} dx_1\cdots dx_{n-1} \nonumber \\
&=\omega_{n-1}\int_U \sqrt{\det(B^T \cdot B)}dx_1\cdots dx_{n-1}  \nonumber \\
&=\omega_{n-1} \lim_{\epsilon\downarrow 0}\frac{\mu_n(Tub_\epsilon W)}{2\epsilon}. \nonumber
\end{align}
Here we have used the fact in the last identity derived in \cite[page 83]{de2013volume}. Finally, the smoothness of $\Sigma$ implies that $\Sigma$ can be covered by such summable open subsets $W$. By the standard argument of partitions of unity, we conclude that
$$\mathcal{V}_{n-1}(\Sigma)=\frac{1}{\omega_{n-1}}\int_{p\in \mathbb{S}^{n-1}}\mathcal{V}_{n-1}(\Sigma|n^\perp(p))dS_p=\lim_{\epsilon\downarrow 0}\frac{\mu_n(Tub_\epsilon \Sigma)}{2\epsilon},$$
and complete the proof.
\end{proof}

\section{generalization to higher codimensions}\label{content2}
In this section, we will project the compact hypersurface $\Sigma$ in $\mathbb{R}^n$ to the lower dimensional subspaces and consider their projected volumes. To be precise, we plan to define the
projected $r$-dimensional volume (projected $r$-volume, in short) of $\Sigma$, $1\leq r\leq n-1$, onto any $r$-dimensional subspace in $\mathbb{R}^n$, and find the recursive formula (see \eqref{recursive1}) for the average projected $r$-volume.
By using the same notation as Santal\'{o}'s in \cite{sors2004integral}, for any point $q\in \mathbb{R}^n$, we denote $L_{r[q]}$ by the $r$-dimensional plane ($r$-plane, in short) in $\mathbb{R}^n$ through $q$, and $L^\perp_{r[q]}$ by its orthogonal complement. Notice that for any point $q\in L_{r[0]}$, there exists the unique affine $(n-r)$-plane $L^\perp_{n-r[q]}$ through $q$ and perpendicular to $L_{r[0]}$.
Indeed, the uniqueness and existence of the affine orthogonal complement can be obtained by parallel shifting $L^\perp_{n-r[0]}$ to $q$, namely, $L^\perp_{n-r[q]}=\{q+v, \text{for any }v\in L^\perp_{n-r[0]}\}$.

Given a $r$-plane $L_{r[0]}$ and a fixed point $q\in \Sigma$. Denote $u$ by the orthogonal projection of $q$ onto $L_{r[0]}$.
Then there exists the unique $(n-r)$-plane containing $p$ and $q$, and perpendicular to $L_{r[0]}$.
Let us denote the $(n-r)$-plane by $L_{n-r[q]}^\perp$. Clearly, $q\in L_{r[0]}\cap L^\perp_{n-r[q]}$.
Now we discuss the dimension of the intersection $T_p\Sigma\cap L^\perp_{n-r[q]}$.
Since the following discussion also holds for any vector subspaces, we may simplify the notations by setting
$$M=T_p\Sigma, \ V=L_{r[0]},\  V^\perp=L^\perp_{n-r[q]}.$$
It is clear that $\displaystyle{\dim(M\cap V^\perp)\leq \min \{\dim(M), \dim(V^\perp)\}}=\dim(V^\perp)=n-r$. Moreover, since the dimension of the sum of $M$ and $V^\perp$, $M+V^\perp=\{m+v, m\in M, v\in V^\perp \}$, is at most $n$, we have
\begin{align*}
n-r-1\leq \dim(M )+\dim(V^\perp)-\dim(M+V^\perp)=\dim(M\cap V^\perp)\leq n-r.
\end{align*}
Thus, there are two cases to be concerned for the dimension $\dim(M\cap V^\perp)$: $n-r-1$ and $n-r$. We will construct the projected $r$-volume forms of $M$ onto $V$ for the first case and show that the all $r$-forms vanish in the second case so that it is measure zero when considering the integral over such points.

\textit{Case 1}. When $\dim(M\cap V^\perp)=n-r-1$, let $\nu$ be the unit normal to $M$. We may choose the orthonormal basis $\{e_\alpha, e, e_\beta\}$ in the space $M\bigoplus span\{\nu\}=\mathbb{R}^n$ and orthonormal basis $\{e_\alpha, u_\delta\}$ in $M$ satisfying the following conditions: ($1\leq \alpha \leq n-r-1, n-r+1\leq \beta\leq n$, and $n-r\leq \delta \leq n-1$):
\begin{enumerate}
  \item $span \{e_\alpha\}=M\cap V^\perp$,
  \item $span \{e_\alpha, e\}=V^\perp$,
  \item $span \{e_\beta\}=V$,
  \item $span \{e_\alpha, u_\delta\}=M$.
\end{enumerate}
Notice that since $\dim(V^\perp \setminus (M\cap V^\perp))=1$, the unit vector $e$ is uniquely determined (up to a sign), independent of the choice of the vectors $e_\alpha, e_\beta$, and $u_\delta$. Since the point $p\in M$, the infinitesimal change $dp$ is still contained in $M$, so we have the vector-valued one-form
$$dp=\sum_{\alpha=1}^{n-r-1} A^\alpha e_\alpha+\sum_{\delta=n-r}^{n-1} B^\delta u_\delta$$ for some connection $1$-forms $A^\alpha, B^\delta$. In addition, $dp$ is a vector in $\mathbb{R}^n$, so it can be written in terms of the linear combination of the basis $\{e_\alpha, e, e_\beta\}$,
namely,
$$dp=\left(\sum_{\alpha=1}^{n-r-1}\omega^\alpha e_\alpha \right)\wedge \omega e \wedge \left(\sum_{\beta=n-r+1}^n \omega^\beta e_\beta\right),$$
for some $1$-forms $\omega^\alpha, \omega$, and $\omega^\beta$.
Thus, for any $\beta=n-r+1, \cdots, n$, one has that
$$\omega^\beta=dp\cdot e_\beta=\sum_{\delta=n-r}^{n-1} B^\delta (u_\delta \cdot e_\beta),$$
and so
\begin{align}\label{areaelement1}
\displaystyle{\bigwedge_{\beta=n-r+1}^{n} \omega^\beta=\Delta\cdot  \bigwedge_{\delta=n-r}^{n-1} B^\delta},
\end{align}
where $\Delta=\det(u_\delta\cdot e_\beta)$ is the determinant of the $(r\times r)$-matrix with entries $u_\delta\cdot e_\beta$.
Notice that the value $\Delta$ satisfies $-1\leq \Delta \leq 1$ and it measures the cosine of the angle between $M\setminus V^\perp$ and $V$, equivalently, the angle between $\nu$ and $e$. Indeed, recall that the Hodge star operator $*$ satisfies the property
$$*(\eta\wedge (*\zeta))=\langle \eta,\zeta \rangle,$$
for any exterior $r$-forms $\eta, \zeta$, where $\langle \  , \ \rangle$ is the inner product for $r$-forms.
Substituting $\eta=\bigwedge_\delta u_\delta$ and $\zeta=\bigwedge_\beta e_\beta$ into the identity, and using the orthogonal decomposition
$$\mathbb{R}^n= M\bigoplus span\{\nu\} =span\{e_\alpha\}\bigoplus span\{u_\delta\}\bigoplus span\{\nu\}\ni e,$$
a straightforward computation shows that
\begin{equation}\label{det2}
\begin{aligned}
\Delta&=\det(u_\delta\cdot e_\beta)\\
&=\langle \bigwedge_\delta u_\delta, \bigwedge_\beta e_\beta \rangle \\
&=*\Big(\bigwedge_\delta u_\delta \wedge (*\bigwedge_\beta e_\beta)\Big)\\
&=*(\bigwedge_\delta u_\delta\wedge e \bigwedge_\alpha e_\alpha)\\
&=(e\cdot  \nu)*(\bigwedge_\delta u_\delta \wedge \nu \bigwedge_\alpha e_\alpha)\\
&=e\cdot \nu,
\end{aligned}
\end{equation}
as desired. To our purpose (see \eqref{defproj}), we will take the absolute value $|\Delta|$ of $\Delta$ in \eqref{areaelement1} such that the $r$-volume form $\wedge_\beta \omega^\beta$ is a positive measure. For more details about the angles between two subspaces with arbitrary dimensions, we refer the reader to  \cite[Theorem 1]{jiang1996angles} and \cite{MIAO199281}.

The geometric meaning of $\eqref{areaelement1}$ can be interpreted as follows:
by identifying the origin $O\in V$ (and so $q$) and $p$ via the parallel transport in $\mathbb{R}^n$, $\bigwedge_\beta \omega^\beta$ (resp. $\bigwedge_{\delta} B^\delta$) is the $r$-dimensional volume element in $V$ at $q$
(resp. in $M\setminus  V^\perp$ at $p$, the subspace in $M$ that is \textit{not} perpendicular to $V$).
The formula \eqref{areaelement1} describes that the projected $r$-form $\bigwedge_\beta \omega^\beta$ is the orthogonal projection of $\bigwedge_\delta B^\delta$ onto the plane $V$, and the projection is independent of the choice of the vectors $e_\alpha, e, e_\beta$, and $u_\delta$.
As a consequence, we have constructed a natural projected $r$-volume element of $M$ onto $V$, and finish the discussion for the first case.

Before advancing to the second case, let us implement the previous construction to the compact hypersurface in $\mathbb{R}^n$.

\begin{prop}\label{prop1}
For any $r$-plane $L_{r[0]}$ through the origin and any point $q\in L_{r[0]}$, there exists the unique $(n-r)$-plane $L_{n-r[q]}^\perp$ through $q$ satisfying
\begin{enumerate}
\item the orthogonal decomposition $L_{r[0]}\bigoplus L_{n-r[q]}^\perp=\mathbb{R}^n$, and
\item\label{prop1vector} if $L_{n-r[q]}^\perp\cap \Sigma \neq \emptyset$, then for any point $p\in L_{n-r[q]}^\perp\cap \Sigma$ with $L_{n-r[q]}^\perp\nsubseteq T_p\Sigma$, there exists the unique (up to a sign) unit vector $e_{p,q}$ in $L_{n-r[q]}^\perp\setminus T_p\Sigma$.
\end{enumerate}
\end{prop}

\begin{proof}
\begin{enumerate}
\item It is clear by the assumption.
\item Suppose $p\in L_{n-r[q]}^\perp\cap \Sigma$. Since $L_{n-r[q]}^\perp\nsubseteq T_p\Sigma$, $\dim(L_{n-r[q]}^\perp\cap T_p\Sigma)=n-r-1$.
Then at $p$ we have the orthogonal decomposition $L_{n-r[q]}^\perp=(L_{n-r[q]}^\perp\cap T_p\Sigma) \bigoplus U$ for some orthogonal complement $U$ of $(L_{n-r[q]}^\perp\cap T_p\Sigma)$ in $L_{n-r[q]}^\perp$.
Besides, since $\dim(U)=\dim(L_{n-r[q]}^\perp)-\dim(L_{n-r[q]}^\perp\cap T_p\Sigma)=1$, we may choose the unique vector starting from the point $p$ with length one (up to a sign) in the $1$-dimensional affine subspace $U$, and the result follows.
\end{enumerate}
\end{proof}

We point out that in \eqref{prop1vector} of Proposition \ref{prop1}, in general, the cross-section $L_{n-r[q]}^\perp\cap \Sigma$ may be comprised of infinitely or finitely many connected components. But for our purpose, we only consider the hypersurface $\Sigma$ with finitely many cross-sections for all $(n-r)$-planes through any point on $\Sigma$.

\begin{prop}\label{prop2}
For any $p\in \Sigma$ and any $(n-r)$-plane $U \nsubseteq T_p\Sigma$ through $p$, there exist the unique $r$-plane $L_{r[0]}$ through the origin and the unique point $q\in L_{r[0]}$  such that
\begin{enumerate}
\item $q\in U\cap L_{r[0]}$ and $U\bigoplus L_{r[0]}=\mathbb{R}^n$ (thus, we may denote $U=L_{n-r[q]}^\perp$ as shown in Proposition \ref{prop1}), and
\item\label{dimencondition1}  $\dim(U\setminus  T_p\Sigma)=1$.
\end{enumerate}
Moreover, by \eqref{dimencondition1}, there exists the unique unit vector $e\in U\setminus T_p\Sigma $ coincided with the vector $e_{p,q}$ constructed in \eqref{prop1vector} of Proposition \ref{prop1}.
\end{prop}

\begin{proof}
\begin{enumerate}
\item The $r$-plane $L_{r[0]}$ can be uniquely obtained by the orthogonal affine subspace to $U$, namely, $L_{r[0]}=\{x\in \mathbb{R}^n, x\cdot y=0 \text{ for all }y\in U \}$, and so the point $q$ is given by the unique point at $L_{r[0]}\cap U$.
\item The transversal assumption $U \nsubseteq T_p\Sigma$ implies the result immediately.
\end{enumerate}
Finally, by \eqref{dimencondition1} we may have the unique (up to a sign) unit vector $e \in U\setminus T_p\Sigma$.
Since the vector $e$ is uniquely determined by $U$,  $L_{r[0]}$, $p$, and $q$, by setting $U=L_{n-r[q]}^\perp$ in Proposition \ref{prop1}, the unique vector $e_{p, q}$ is exactly same as the vector $e$.
\end{proof}

\begin{remark}
For any fixed point $p\in \mathbb{R}^n$, the natural orthogonal decomposition for $\mathbb{R}^n$ implies that there exists a one-to-one correspondence that assigns the $r$-plane $L_{r[0]}$ a $(n-r)$-plane $L^\perp_{n-r[q]}$ through $p$, where $q$ is the orthogonal projection of $p$ onto $L_{r[0]}$.
Let $G_{n,r}$ be the Grassmannian, the set of all $r$-subspaces in $\mathbb{R}^n$. According to Proposition \ref{prop2}, if $p\in \Sigma$, the map
\begin{align}\label{rcl1}
\begin{array}{rcl}
\phi_p: G_{n,r}  & \rightarrow  &  G_{n-r, 1} \\
L_{r[0]} & \mapsto      &  \phi(L_{r[0]})=e\in L_{n-r[q]}^\perp \setminus T_p\Sigma
\end{array}
\end{align}
is a bijection except for the $r$-planes $L_{r[0]}$ perpendicular to $T_p\Sigma$ (namely, $L^\perp_{n-r[q]}\subseteq T_p\Sigma$ for some point $q\in L_{r[0]}$ in Proposition \ref{prop1} (\ref{prop1vector})).
% Indeed, suppose $e$ and $e'$ are the unit vectors at $p\in \Sigma$ obtained from the $r$-planes $L_{r[0]}$ and $L_{r[0]}^\prime$, respectively, by Proposition \ref{prop2} with $e=e'$. Since $L_{r[0]}$ and $L_{r[0]}^\perp$ both are not parallel to $T_p\Sigma$, their corresponding orthogonal complements $L_{n-r[q]}^\perp$ and $L_{n-r[q']}^\perp$ are not perpendicular to $T_p\Sigma$ for some $q\in L_{r[0]}$, $q'\in L_{r[0]}^\perp$. If $q=q'$, then $L_{r[0]}=L_{r[0]}^\prime$ and the proof is done. If $q\neq q'$, then there exists a nonzero angle between those $(n-r)$-planes. Thus, $\dim(L_{n-r[q]}^\perp \cap L_{n-r[q']}^\perp)=n-r-1=\dim$ $e, e'$ are in the orthogonal complement of $L_{n-r[q]}$ and $L_{n-r[q']}$ respectively, and $T_p\Sigma$
\end{remark}

Next, let us continue to discuss the second case of $\dim(M\cap V^\perp)$.

\textit{Case 2}. When $\dim(M\cap V^\perp)=n-r$, it means that $V^\perp$ is contained in $M$. We claim that, by the similar construction in \textit{Case 1}, the $r$-form $\wedge_\beta \omega^\beta$ vanishes.
%Since $\nu \perp M$ (and so $\nu \perp V^\perp $), $\nu \in (V^\perp)^\perp= V$. We may choose one of the vectors $\{ e_\alpha\}_{\alpha=1}^{n-r}$ to be the unit normal vector $\nu$, say $\nu=e_{n-r}$.
One may choose the orthonormal basis $\{e_\alpha, e_\beta\}$ in $M\bigoplus span\{\nu\}= \mathbb{R}^n$ and $\{ e_\alpha, u_\delta\}$ in $M$ satisfying ($1\leq \alpha \leq n-r$, $n-r+1\leq \beta \leq n$, $n-r+1\leq \delta\leq  n-1$)
\begin{enumerate}
\item $span\{e_\alpha\}= M\cap V^\perp = V^\perp$,
\item $span\{e_\beta\}=V$,
\item $span\{e_\alpha, u_\delta\}=M$.
\end{enumerate}
Similar to \textit{Case 1}, on one hand, since $p=M\cap V^\perp$, $dp=\sum_{i=1}^n\omega ^i e_i$ for some connection $1$-forms $\omega^i$. On the other hand, by writing
$$dp=\sum_{\alpha=1}^{n-r} A^\alpha e_\alpha + \sum_{\delta=n-r+1}^{n-1} B^\delta u_\delta,$$ and taking the inner product with $e_\beta$, $n-r+1\leq \beta \leq  n$, one has
\begin{align*}
\omega^\beta =dp\cdot e_\beta=\sum_{\delta=n-r+1}^{n-1} B^\delta (u_\delta \cdot e_\beta).
\end{align*}
We deduce
\begin{align}\label{areaelement2}
\bigwedge_{\beta=n-r+1}^{n} \omega^\beta=\bigwedge_{\beta=n-r+1}^{n} \Big( \sum_{\delta=n-r+1}^{n-1} B^\delta (u_\delta \cdot e_\beta)\Big)=0.
\end{align}
The last equality holds since the wedge product makes a $r$-form from $(r-1)$ one-forms $B^\delta$, and there must be some $B^\delta$'s repeated. This finishes the discussion for \textit{Case 2}.

In contrast to \textit{Case 1}, \textit{Case 2} shows that if $V^\perp \subset M$ (equivalently, $V$ is perpendicular to $M$), the orthogonal contribution of any $r$-volume element in $M$ onto $V$ is zero.
As a consequence of both cases, when considering the integral over all projected $r$-volumes $\wedge_{\beta=n-r+1}^{n}\omega^\beta$, \eqref{areaelement1} and \eqref{areaelement2} suggest that
we may ignore \textit{Case 2}, and only consider  \textit{Case 1}, $\dim(M\cap V^\perp)=n-r-1$.

According to the discussion above, we give a definition for the weighted projected area of any compact hypersurface $\Sigma$ onto any subspace $L_{r[0]}$ of lower dimension.

\begin{defn}\label{codimprojected}
Given a compact hypersurface $\Sigma$ in $\mathbb{R}^n$ and any $r$-plane $L_{r[0]}$ through the origin, $1\leq r\leq n-1$. The (orthogonal) weighted projected $r$-volume $\mathcal{V}_r(\Sigma|L_{r[0]})$ of $\Sigma$ onto $L_{r[0]}$ is defined by
\begin{align}\label{defproj}
\mathcal{V}_r(\Sigma| L_{r[0]})
%=\int_{p\in \Sigma} L_{pq^{-1}}^* \Big( \bigwedge_{j=n-r+1}^n \omega^j(p) \Big)
=\frac{1}{2}\int_{p\in \Sigma} L_{pq^{-1}}^* \Big( |\Delta(p)| \bigwedge_{\delta=1}^r B^\delta(p)\Big),
\end{align}
where $L_{pq^{-1}}^*$ is the pullback of the left translation $L_{pq^{-1}}$, $q$ is the orthogonal projection of $p\in \Sigma$ onto $L_{r[0]}$,
$\Delta(p)$ is the angle between the unit normal $\nu$ to $T_p\Sigma$ and the unique vector $e$ in $L^\perp_{n-r[q]}$ defined in \eqref{rcl1},
$\bigwedge_{\delta=1}^r B^{\delta}(p)$ is the $r$-volume in $T_p\Sigma\setminus L_{n-r[q]}^\perp$.
Also, the mean value of the projected $r$-volumes $\mathcal{V}_r(\Sigma|L_{r[0]})$ is defined by
\begin{align}\label{defmean}
E(\mathcal{V}_r(\Sigma|L_{r[0]}))=\frac{\int_{L_{r[0]}\in G_{n,r}} \mathcal{V}_r(\Sigma| L_{r[0]}) dL_{r[0]}}{m(G_{n,r})},
\end{align}
where $G_{n,r}$ is the Grassmannian consisting of all $r$-dimensional subspaces in $\mathbb{R}^n$, $dL_{r[0]}$ is the invariant density defined below \eqref{invdensity}, and $m(G_{n,r})$ is the volume of the Grassmannian $G_{n,r}$.
\end{defn}

We point out that when $r=n-1$, $T_pM\cap L_{1[q]}^\perp=\{p\}$, and so $\wedge_{\delta=1}^{n-1} B^\delta$ is the surface area element of $\Sigma$ at $p$. Thus, \eqref{defproj} coincides with \eqref{gendefn1}.
When $r=1$, we have the mean width \eqref{defmean} for arbitrary $\Sigma$.

Recall \cite[page 202]{sors2004integral} that the invariant density of $L_{r[0]}$ is given by
\begin{align}\label{invdensity}
\displaystyle{dL_{r[0]}=\bigwedge_{\genfrac{}{}{0pt}{3}{1\leq h\leq n-r}{n-r+1 \leq \beta\leq n}} \omega^h_\beta},
\end{align}
where
\begin{align*}
\omega_\beta^h=
\left\{
\begin{array}{ll}
de_\beta \cdot e_\alpha & \text{ if } 1\leq h\leq n-r-1,\\
de_\beta \cdot e        & \text{ if } h=n-r,
\end{array}
\right.
\end{align*}
(use the same indices $\alpha, \beta$ as in \textit{Case 1}) and the identity (12.36) in \cite[page 203]{sors2004integral} gives
\begin{align}\label{intevalue}
\int_{G_{r,q}}dL_{r[q]}=\frac{O_{n-q-1}O_{n-q-2}\cdots O_{n-r}}{O_{r-q-1}O_{r-q-r}\cdots O_1O_0}
\end{align}
for $0\leq q<r\leq n-1$.
We also point out that the volume $m(G_{n,r})=\frac{O_{n-1}\cdots O_{n-r}}{O_{r-1}\cdots O_1 O_0}$, where $O_r=\frac{2\pi^{(r+1)/2}}{\Gamma((r+1)/2)}$ is the surface area of the unit ball in $\mathbb{R}^{r+1}$
and $\Gamma$ denotes the gamma function. For instance, $O_0=2$ (by convention), $O_1=2\pi$ and $O_2=4\pi$. Notice that the $r$-form $\bigwedge_{\beta=n-r+1}^n \omega^\beta$ (resp. $\bigwedge_{\delta=1}^r B^\delta$) is the $r$-dimensional
volume of the infinitesimal parallelotope in $L_{r[q]}$ (resp. in $T_p\Sigma\setminus L_{{n-r}[q]}^\perp$).
Moreover, if $\Sigma$ is convex, then the definitions \eqref{defproj} and \eqref{defmean} coincide with (13.1) and (13.2) respectively in \cite{sors2004integral}, and so we have had the generalized projected $r$-volumes for arbitrary hypersurface $\Sigma$.

The rest of this paper will be devoted to the derivation of a recursive formula for the integral of projected $r$-volumes.
Recall that \cite[page 216 (13.2)]{sors2004integral} the integral of the projected $(n-r)$-volume of a convex body $K$ (called the mean $(n-r)$-volume in short) is defined by
\begin{align}\label{meanvolume1}
I_r(K)=\int_{G_{n,r}} V(K'_{n-r}) dL_{r[0]}=\int_{G_{n,n-r}} V(K'_{n-r}) dL_{n-r[0]},
\end{align}
where $K'_{n-r}$ is the convex set of all intersection points of $L_{n-r[0]}$ with the $r$-planes perpendicular to $L_{n-r[0]}$ through each point of $K$ and $V(K'_{n-r})$ is the $(n-r)$-volume of $K'_{n-r}$.
Using our notation, it means that $V(K'_{n-r})=\mathcal{V}_{n-r}(K|L_{n-r[0]})$.
The definition \eqref{meanvolume1} of $I_r(K)$ can be generalized to arbitrary compact submanifold $K$ with smooth boundary $\Sigma=\partial K$ if the $(n-r)$-volume $V(K'_{n-r})$ is replaced by that of \eqref{defproj}, namely,
the mean $(n-r)$-volume of $\Sigma$ is defined by
\begin{align}\label{meanvolume2}
I_r(\Sigma)=\int_{G_{n,n-r}}\mathcal{V}_{n-r}(\Sigma|L_{n-r[0]}) dL_{n-r[0]}
\end{align}
which is exactly same as shown in \eqref{defmean} (up to a constant $m(G_{n,r})$.

\begin{Thm}\label{mainthm2}
Let $K$ be an $n$-dimensional compact submanifold in $\mathbb{R}^n$ with boundary $\Sigma=\partial K$. Denote $I_r(K)$ by the mean $(n-r)$-volume as defined in \eqref{meanvolume2}. Then we have the recursive formula
\begin{align}\label{recursive1}
I_r(K)=\frac{2}{O_{r-1}}\int_{G_{n,n-1}} I^{(n-1)}_{r-1}(K'_{n-1}) dL_{n-1[0]},
\end{align}
where $O_{r-1}$ is the surface area of the unit ball in $\mathbb{R}^{r}$, and $I^{(n-1)}_{r-1}(K'_{n-1})$ is the mean $(r-1)$-volume of the projection $K'_{n-1}$ of $K$ onto $L_{n-1[0]}$.
\end{Thm}

\begin{remark}
In \cite[page 217]{sors2004integral}, the author derived the same recursive formula (see the identity (13.7) there) with the assumption that $K$ is a convex body in $\mathbb{R}^n$. We observe that the similar argument can be applied even for nonconvex domains $K$ when the new concept for projected $r$-volumes (i.e. \eqref{gendefn1} and \eqref{defproj}) is introduced. The main idea of the proof of Theorem \ref{mainthm2} is based on the identities \eqref{density1} and \eqref{density2}, which are irrelevant to the convexity for $K$.
\end{remark}

\begin{proof}
Given a $r$-plane $L_r$, $1\leq r\leq n-1$, in $\mathbb{R}^n$. Denote $L_{i+1}^{(r)}$ by the $(i+1)$-plane contained in $L_r$ for $i+1\leq r\leq n-1$. In \cite[page 207 $(12.53)$]{sors2004integral},
the author considered the density for the sets of pairs of linear subspaces $(L_r, L_{i+1}^{(r)})$ and has the identity
\begin{align}\label{density1}
dL_{i+1}^{(r)}\wedge dL_r^* = d L_{r[i+1]}\wedge dL_{i+1},
\end{align}
where $dL_{r}^*$ is the density of the oriented $r$-plane $L_r$ and $dL_{r[i+1]}$ is the density for $r$-planes about a fixed $(i+1)$-plane.
If we consider the linear spaces through the fixed origin $O$ in $\mathbb{R}^n$, \eqref{density1} still holds and may be written
\begin{align}\label{density2}
dL_{i+1[0]}^{(r)}\wedge dL_{r[0]}^* = dL_{r[i+1]}\wedge dL_{i+1[0]}.
\end{align}
In particular, when $L_r$ is a hyperplane and $L_{i+1[0]}^{(r)}$ is of maximal dimension in $L_r$, namely, $r=n-1$, $i+1=r$, \eqref{density2} becomes
\begin{align}\label{density3}
dL_{r[0]}^{(n-1)}\wedge dL_{n-1[0]}^* =dL_{n-1[r]}\wedge dL_{r[0]}.
\end{align}
Similarly, when restrict to the hyperplane $L_{n-1[0]}$ (namely, substitute $n$ by $n-1$, and $i+1$ by $r-1$) in \eqref{density2}, one has
\begin{align}\label{density4}
dL_{r-1[0]}^{(r)} \wedge dL_{r[0]}^{*(n-1)} = dL_{r[r-1]}^{(n-1)} \wedge dL_{r-1[0]}^{(n-1)},
\end{align}
where the superscripts ${(n-1)}$ emphasize that the sub-planes considered here are contained in the plane $L_{n-1[0]}$. For instance, $ dL_{r[0]}^{*(n-1)}$ is the density of oriented $r$-plane through the origin $O$ contained in $L_{n-1[0]}$.
Multiplying \eqref{density3} by $dL^r_{r-1[0]}$, \eqref{density4} by $dL_{n-1[0]}$, and using the fact that an oriented plane is equivalent to two unoriented planes such that one has $dL_{n-1[0]}^*=2dL_{n-1[0]}$ and $dL^{*(n-1)}_{r[0]}=2dL^{(n-1)}_{r[0]}$,
we have reached
\begin{align}\label{density6}
dL_{r[r-1]}^{(n-1)}\wedge dL_{r-1[0]}^{(n-1)}\wedge dL_{n-1[0]} = dL_{r-1[0]}^{(r)} \wedge  dL_{n-1[r]}\wedge  dL_{r[0]}.
\end{align}

Now let us integrate over all the pairs of $L_{r-1[0]}$ and $L^{(n-1)}_{r[0]}$. Let $\Sigma_{n-1}^\prime$ be the projected $(n-1)$-volume of $\Sigma$ onto $L_{n-1[0]}$.
Notice that the projected $(n-r)$-volume of $\Sigma$ onto $L_{n-r[0]}$ is equal to the projected $(n-r)$-volume of $\Sigma_{n-1}^\prime$ onto $L_{n-r[0]}$ (counted for multiplicities). On the one hand, the integral of the left-hand side of \eqref{density6} becomes
\begin{align}\label{rhs}
&\int_{G_{n,r}}\int_{G_{r,r-1}}\int_{G_{n-1,r}}  V(K'_{n-r}) dL_{r[r-1]}^{(n-1)}\wedge dL_{r-1[0]}^{(n-1)}\wedge dL_{n-1[0]}   \\
&=\int_{G_{n,r}}I^{(n-1)}_{r-1}(K'_{n-1})  dL_{n-1[0]}\wedge dL_{r[r-1]}^{(n-1)} \nonumber \\
&=\int_{G_{r,r-1}^{(n-1)}}\int_{G_{n,n-1}} I^{(n-1)}_{r-1}(K'_{n-1})dL_{n-1[0]}\wedge dL_{r[r-1]}^{(n-1)}\nonumber \\
&=\int_{G_{n,n-1}}I^{(n-1)}_{r-1}(K'_{n-1}) dL_{n-1[0]}\cdot \int_{G_{r,r-1}^{(n-1)}} dL^{(n-1)}_{r[r-1]} \nonumber \\
&=\frac{O_{n-r-1}}{2}\int_{G_{n,n-1}}I^{(n-1)}_{r-1}(K'_{n-1})dL_{n-1[0]}, \nonumber
\end{align}
here we have used that $\int_{G_{r,r-1}}\int_{G_{n-1,r}}=\int_{G_{n-1,r-1}}$, $\int_{G_{n,r}}=\int_{G_{n,n-1}}\int_{G^{(n-1)}_{r,r-1}}$ in the first two identities, and \eqref{intevalue} in the last identity $\frac{O_{n-r-1}}{2}=\int_{G_{r,r-1}^{(n-1)}} dL^{(n-1)}_{r[r-1]}$.
On the other hand, we deduce the integral of the right-hand side of \eqref{density6}
\begin{align}\label{lhs}
\indent &\int_{G_{n,r}}\int_{G_{r,r-1}}\int_{G_{n-1,r}} V(K'_{n-r})   dL_{r-1[0]}^{(r)} \wedge  dL_{n-1[r]}\wedge  dL_{r[0]}\\
&=\int_{G_{r,r-1}}\int_{G_{n-1,r}} I_r(K) dL_{r-1[0]}^{(r)} \wedge  dL_{n-1[r]} \nonumber \\
&=\frac{O_{r-1}}{2}\frac{O_{n-r-1}}{2}I_r(K). \nonumber
\end{align}
Again we have used \eqref{intevalue} to have $\int_{G_{r,r-1}}dL_{r-1[0]}^{(r)}=\frac{O_{r-1}}{O_0}$ and $\int_{G_{n-1,r}} dL_{n-1[r]}= \frac{O_{n-r-1}}{O_0}$ in the last equality. Combining \eqref{rhs} and \eqref{lhs} to have the recursive formula \eqref{recursive1}.
\end{proof}

\bibliography{mybib1}

\begin{bibdiv}
\begin{biblist}

\bib{aminov2001geometry}{book}{
      author={Aminov, Yu},
       title={The geometry of submanifolds},
   publisher={CRC Press},
        date={2001},
}

\bib{bouafia2015integral}{article}{
      author={Bouafia, Philippe},
      author={De~Pauw, Thierry},
       title={Integral geometric measure in separable banach space},
        date={2015},
     journal={Mathematische Annalen},
      volume={363},
      number={1},
       pages={269\ndash 304},
}

\bib{chakerian1966sets}{article}{
      author={Chakerian, Gulbank},
       title={Sets of constant width},
        date={1966},
     journal={Pacific Journal of Mathematics},
      volume={19},
      number={1},
       pages={13\ndash 21},
}

\bib{de2013volume}{article}{
      author={de~Jong, Theo},
       title={Volume of submanifolds},
        date={2013},
     journal={Mathematische Semesterberichte},
      volume={60},
      number={1},
       pages={81\ndash 83},
}

\bib{favard1932definition}{article}{
      author={Favard, J},
       title={Une definition de la longueur et de l'aire},
        date={1932},
     journal={CR Acad. Sci. Paris},
      volume={194},
       pages={344\ndash 346},
}

\bib{gardner1995geometric}{book}{
      author={Gardner, Richard~J},
       title={Geometric tomography},
     edition={2},
   publisher={Cambridge University Press Cambridge},
        date={2006},
      volume={6},
}

\bib{huang2021cauchy}{article}{
      author={Huang, Yen-Chang},
       title={Cauchy’s surface area formula in the heisenberg groups},
        date={2021},
     journal={Revista Matematica Iberoamericana},
}

\bib{jiang1996angles}{article}{
      author={Jiang, Sheng},
       title={Angles between euclidean subspaces},
        date={1996},
     journal={Geometriae Dedicata},
      volume={63},
      number={2},
       pages={113\ndash 121},
}

\bib{klain1997}{book}{
      author={Klain, Daniel~A},
      author={Rota, Gian-Carlo},
      author={others},
       title={Introduction to geometric probability},
   publisher={Cambridge University Press},
        date={1997},
}

\bib{MIAO199281}{article}{
      author={Miao, Jianming},
      author={Ben-Israel, Adi},
       title={On principal angles between subspaces in rn},
        date={1992},
        ISSN={0024-3795},
     journal={Linear Algebra and its Applications},
      volume={171},
       pages={81\ndash 98},
  url={https://www.sciencedirect.com/science/article/pii/0024379592902515},
}

\bib{sors2004integral}{book}{
      author={Santal{\'o}, Luis~Antonio},
       title={Integral geometry and geometric probability},
   publisher={Cambridge university press},
        date={2004},
}

\bib{schneider2014convex}{book}{
      author={Schneider, Rolf},
       title={Convex bodies: the brunn--minkowski theory},
   publisher={Cambridge university press},
        date={2014},
      number={151},
}

\bib{tsukerman2017brunn}{article}{
      author={Tsukerman, Emmanuel},
      author={Veomett, Ellen},
       title={Brunn-minkowski theory and cauchy's surface area formula},
        date={2017},
     journal={The American Mathematical Monthly},
      volume={124},
      number={10},
       pages={922\ndash 929},
}

\bib{zeng2016new}{article}{
      author={Zeng, Chunna},
      author={Bai, Shikun},
      author={Tong, Yin},
       title={A new integral formula for the angle between intersected
  submanifolds},
        date={2016},
     journal={Journal of Inequalities and Applications},
      volume={2016},
      number={1},
       pages={1\ndash 9},
}

\bib{zhang1999dual}{article}{
      author={Zhang, Gaoyong},
       title={Dual kinematic formulas},
        date={1999},
     journal={Transactions of the American Mathematical Society},
      volume={351},
      number={3},
       pages={985\ndash 995},
}

\end{biblist}
\end{bibdiv}
\end{document}